\documentclass[10pt]{amsart}
\usepackage{url,graphicx}
\usepackage[dvips]{epsfig}
\usepackage{latexsym,color}
\usepackage{amscd,amssymb,verbatim}
\newcommand{\be}{\begin{enumerate}}
\newcommand{\ee}{\end{enumerate}}

\newcommand{\cp}{{\mathcal P}}

\newcommand{\rn}{{\mathbb N}}
\newcommand{\rr}{{\mathbb R}}

\newcommand{\cz}{{\mathbb Z}}
\newcommand{\da}{\Delta}

\newcommand{\ra}{\rightarrow}
\newcommand{\nin}{\noindent}
\newcommand{\pr}{\noindent{\bf Proof. }}
\newcommand{\sm}{\setminus}

\newcommand{\bo}{\partial}
\newcommand{\gt}{\rightarrow}
\newcommand{\gti}{\gt\infty}
\newcommand{\lni}{\lim_{n\gti}}
\newcommand{\prob}{\textrm{Prob}\,}
\newcommand{\npart}{\textrm{part}\,}
\newcommand{\ynpd}{Y_{n,p,d}}
\newcommand{\ynpone}{Y_{n,p,1}}
\newcommand{\ynptwo}{Y_{n,p,2}}
\newcommand{\supp}{\textrm{supp}\,}
\newcommand{\zz}{{\mathbb Z}_2}

\newcommand{\ti}{\tilde}
\newcommand{\wti}{\widetilde}

\newtheorem{thm}{Theorem}[section]
\newtheorem{df}[thm]{Definition}
\newtheorem{lm}[thm]{Lemma}

\numberwithin{equation}{section}
\numberwithin{figure}{section}

\begin{document}

\title[The threshold for vanishing of the top homology] {The threshold
  function for vanishing of the top homology group of random
  $d$-complexes}

\author{Dmitry N. Kozlov}
\address{Department of Mathematics, University of Bremen, 28334 Bremen,
Federal Republic of Germany}
\email{dfk@math.uni-bremen.de}
\thanks {This research was supported by University of Bremen, 
as part of AG CALTOP}
\thanks {Note that the title and MSC numbers are provisional}

\date\today

\begin{abstract}
For positive integers $n$ and $d$, and the probability function $0\leq
p(n)\leq 1$, we let $\ynpd$ denote the probability space of all at
most $d$-dimensional simplicial complexes on $n$ vertices, which
contain the full $(d-1)$-dimensional skeleton, and whose $d$-simplices
appear with probability $p(n)$. In this paper we determine the
threshold function for vanishing of the top homology group in $\ynpd$,
for all $d\geq 1$.
\end{abstract}

\maketitle

\section{Thresholds for vanishing of the $(d-1)$st homology group
of random $d$-complexes}

In 1959 Erd\H{o}s and R\'enyi have defined a~natural model for random
graphs which has since become classical. In this model, which we call
$Y_{n,p,1}$, the random graph always has $n$ vertices, where $n$ is
fixed, and the edges are chosen uniformly at random with probability
$p$. One of their main results concerning $Y_{n,p,1}$ was the
discovery of the threshold function for the connectivity of the
graph. More precisely, reformulated in our language, they have shown
the following theorem.

\begin{thm} {\rm(Erd\H{o}s-R\'enyi Theorem, \cite{ER}).}

\nin Assume that $w(n)$ is any function $w:\rn\ra\rr$, such that $\lni
w(n)=\infty$, and $p=p(n)$ is probability depending on~$n$, then we
have

\begin{enumerate}
\item if $p(n)=\frac{\log n-w(n)}{n}$, then
  $\lni\prob(\wti\beta_0(\ynpone;\zz)>0)=1$;
\item if $p(n)=\frac{\log n+w(n)}{n}$, then
  $\lni\prob(\wti\beta_0(\ynpone;\zz)=0)=1$.
\end{enumerate}
\end{thm}

More recently, the two-dimensional analog $Y_{n,p,2}$ of
Erd\H{o}s-R\'enyi model was considered by Linial-Meshulam in~\cite{LM},
and, further, the $d$-dimensional model $\ynpd$, for $d\geq 3$, was
considered by Meshulam-Wallach in~\cite{MW}.

In these generalizations, the graphs are replaced with simplicial
complexes of dimension at most $d$, on $n$ vertices, where all
simplices of dimension $d-1$ or less are required to be in the
complex, and the simplices of dimension $d$ are chosen uniformly at
random with probability~$p$. The combined work of Linial-Meshulam and
Meshulam-Wallach yields threshold functions for the vanishing of the
$(d-1)$th homology group of $\ynpd$ with coefficients in a~finite
abelian group. Specifically, the following is known.

\begin{thm} {\rm{(Linial-Meshulam, \cite{LM}; 
Meshulam-Wallach, \cite{MW}).}}

\nin Assume that $w(n)$ is any function $w:\rn\ra\rr$, such that $\lni
w(n)=\infty$, and $p=p(n)$ is probability depending on~$n$, and $F$ is
a~finite abelian group. Then we have

\begin{enumerate}
\item if $p(n)=\frac{d\log n-w(n)}{n}$, then
  $\lni\prob(H_{d-1}(\ynpd;F)\neq 0)=1$;
\item if $p(n)=\frac{d\log n+w(n)}{n}$, then
  $\lni\prob(H_{d-1}(\ynpd;F)=0)=1$.
\end{enumerate}
\end{thm}

Curiously, the methods of \cite{LM,MW} do not easily extend to the
case of integer coefficients, and finding the threshold functions for
the vanishing of $H_{d-1}(\ynpd;\cz)$ remains open even for the case
$d=2$.

On the other hand, the threshold for vanishing of the fundamental
group of $\da\in\ynptwo$ is well understood due to work of Babson,
Hoffman, and Kahle. The following deep result can be found
in~\cite{BHK}.

\begin{thm} {\rm(Babson, Hoffman, Kahle, \cite[Theorem 1.3]{BHK}).}

\nin If $w(n)$ is a~function, such that $\lni w(n)=\infty$, and
$p(n)\geq \left(\frac{3\log n+w(n)}{n}\right)^{1/2}$, then
$\prob(\pi_1(\ynptwo)=0)=1$.
\end{thm}

Since the simplicial complexes $\da$ in $\ynpd$ have dimension at most
$d$, and are, on the other hand, required to contain full
$(d-1)$-dimensional skeleton, we have $H_i(\da;F)=0$, for all $i\neq
d-1,d$, where $F$ is an~arbitrary abelian group. In this paper we
complement the study undertaken by Linial-Meshulam and
Meshulam-Wallach, by computing the threshold functions for the
vanishing of the top dimensional homology.

\section{Terminology and the formulation of the main result}

We start by recalling some standard notations.  For a~positive integer
$n$, we let $\da_n$ denote the full $(n-1)$-dimensional simplex. Given
a simplicial complex $\da$, and a~nonnegative integer~$d$, we let
$\da^{(d)}$ denote the $d$-dimensional skeleton of $\da$, and we let
$\da(d)$ denote the set of the $d$-simplices of~$\da$. Furthermore,
for an~arbitrary abelian group $F$, we let $B_{d-1}(\da;F)$ denote the
subspace of $C_{d-1}(\da;F)$ generated by the boundaries of the
$d$-simplices from~$\da$, and we let $Z_d(\da;F)$ denote the subspace
of $C_d(\da;F)$ consisting of the cycles. Finally, for a~$d$-chain
$\sigma\in C_d(\da;F)$ we let $\supp\sigma$ denote the subset of
$\da(d)$ consisting of all $d$-simplices appearing with non-zero
coefficients in~$\sigma$. We also assume familiarity with
Bachmann-Landau notations for the asymptotic behavior of functions.

For positive integers $n$ and $d$, and a~real number $0\leq p\leq 1$,
we let $\ynpd$ denote the probability space of all at most
$d$-dimensional simplicial complexes on $n$ vertices, which contain
the full $(d-1)$-dimensional skeleton, and whose $d$-simplices appear
with probability $p$. Formally, the underlying set of $\ynpd$ consists
of all simplicial complexes $\da$, such that
$\da_n^{(d-1)}=\da^{(d-1)}$, and $\da(d+1)=\emptyset$; clearly there
are $2^{n\choose{d+1}}$ of them. The probability associated to each
$\da$ is
$p^{\left|\da(d)\right|}(1-p)^{{n\choose{d+1}}-\left|\da(d)\right|}$. When
the values $n$, $p$, and $d$ are fixed, and $S$ is some set of
simplices of $\da_n$, we shall write $\prob(S)$ to denote the
probability that all of the simplices from $S$ are present in the
simplicial complex sampled from $\ynpd$.

To work with the probability space $\ynpd$ we shall use the following
notations. We write $\da\in\ynpd$ when we sample a~simplicial complex
from $\ynpd$. For any integer $i$, and any field $F$, we write
$\beta_i(\ynpd;F)$ to denote the expectation of the $i$th Betti number
in the probability space $\ynpd$. We also write
$\prob(\beta_i(\ynpd;F)=0)$, and $\prob(\beta_i(\ynpd;F)>0)$ to denote
the probabilities that the $i$th Betti number of $\da\in\ynpd$ is
equal to $0$, correspondingly is strictly larger than~$0$. Similarly,
for an~arbitrary abelian group $F$, we write $\prob(H_i(\ynpd;F)=0)$,
and $\prob(H_i(\ynpd;F)\neq 0)$, to denote the probabilities that the
$i$th homology group of $\da\in\ynpd$ is trivial, correspondingly
nontrivial.

To keep our argument as simple as possible, we shall initially
restrict ourselves to $\zz$-coefficients. The adjustments needed to
handle the general case will follow in Section~\ref{sect:5}.

\begin{thm}\label{thm:main}
The probability $p(n)=\Theta\left(\frac{1}{n}\right)$ is the threshold
probability for vanishing of the top homology of the random simplicial
$d$-complex.  More precisely, assume that $p=p(n)=w(n)/n$, and $d\geq
1$, then we have

\begin{enumerate}
\item if $\lni w(n)=0$, then $\lni\prob(\beta_d(\ynpd;\zz)=0)=1$;
\item if $\lni w(n)=\infty$, then $\lni\prob(\beta_d(\ynpd;\zz)>0)=1$.
\end{enumerate}
\end{thm}

Before proceeding with the proof, we need two more pieces of notation.

\begin{df}\label{df:sigma0}
For an~arbitrary positive integer $d$, let $\Sigma_d$ denote the
subset of $C_{d-1}(\da_n;\zz)\times 2^{\da_n(d)}\times\cz_{\geq 0}$
defined by the following: $(\sigma,S,\lambda)\in\Sigma_d$ if and only
if $|\supp\sigma|>(\lambda-1)(d+1)$.
\end{df}

\noindent
In particular, $(0,S,\lambda)\in\Sigma_d$ implies $\lambda=0$.

\begin{df} \label{df:sigma}
For $(\sigma,S,\lambda)\in\Sigma_d$, we define
$\rho(\sigma,S,\lambda)$ to be the probability that $\da\in\ynpd$
satisfies the following two conditions:
\begin{enumerate}
\item[(1)] $\da$ contains $\sigma$ in its boundary set,
i.e., $\sigma\in B_{d-1}(\da;\zz)$;
\item[(2)] the sets $\da(d)$ and $S$ are disjoint.
\end{enumerate}
\end{df}

\noindent
So, informally speaking, a~collection of the simplices from $\da$ can
be used to complement $\sigma$ to a~$d$-cycle, avoiding the
$d$-simplices from~$S$.

For future reference, we record a few simple properties of
$\rho(-,-,-)$.
\begin{lm} \label{lm:easy} $\,$
\noindent
\begin{enumerate}
\item[(1)] For any $(\sigma,S_1,\lambda),(\sigma,S_2,\lambda)\in\Sigma_d$, 
the set inclusion $S_2\subseteq S_1$ implies
$\rho(\sigma,S_1,\lambda)\leq\rho(\sigma,S_2,\lambda)$;
\item[(2)] we have $\rho(\sigma,S,\lambda_1)=\rho(\sigma,S,\lambda_2)$, 
whenever $(\sigma,S,\lambda_1),(\sigma,S,\lambda_2)\in\Sigma_d$;
\item[(3)] whenever $(\sigma,\da_n(d),\lambda)\in\Sigma_d$, and
  $\sigma\neq 0$, we have $\rho(\sigma,\da_n(d),\lambda)=0$;
\item[(4)] we have $\rho(\sigma,S,\lambda)=0$, for all 
$(\sigma,S,\lambda)\in\Sigma_d$, such that $\bo\sigma\neq 0$.
\end{enumerate}
\end{lm}
\pr The first condition holds simply because in $\ynpd$ it is less
probable that a~simplicial complex satisfies a~(possibly) more
stringent set of conditions. The second condition is
straightforward. The third condition is true since in this case
$\da(d)$ must be empty. Finally, the fourth condition holds since
the square of the differential in a~chain complex is equal
to~$0$. \qed


\section{Proof of the first part of Theorem~\ref{thm:main}}

We start with the first part of Theorem~\ref{thm:main}, which is more
difficult (and more interesting). Its proof relies on the following
lemma, which might also be useful in its own right.

\begin{lm} \label{star}
Let us fix positive integers $n$ and $d$, and a~probability $1\geq
p\geq 0$, such that $d\geq 2$, $n\geq d+1$, and $pn<1$. Set
$w:=pn$. For any $(\sigma,S,\lambda)\in\Sigma_d$ we have
\begin{equation} \label{eq:est}
\rho(\sigma,S,\lambda)\leq 
c(d,\lambda)\, p^\lambda/(1-w)^\lambda,
\end{equation}
where $c(d,\lambda)=(d+1)^\lambda\lambda!$.
\end{lm}

The case $S=\emptyset$ is of special interest to us and we adopt the
abbreviated notation
$\rho(\sigma,\lambda):=\rho(\sigma,\emptyset,\lambda)$.

\vspace{5pt}

\nin {\bf Proof of Lemma~\ref{star}.} 
By Lemma~\ref{lm:easy}(4) we can always assume that $\bo\sigma=0$, as
otherwise the left hand side of~\eqref{eq:est} is equal to~$0$.

We shall use induction on $\lambda$. The base of induction is
$\lambda=0$. In this case $c(d,0)=1$ for all~$d$, and the right hand
side of \eqref{eq:est} is equal to $1$; hence the inequality is
trivially satisfied.

To prove the induction step, let us now assume that $\lambda\geq 1$,
and that the inequality~\eqref{eq:est} has been shown for all
$\ti\lambda$, such that $0\leq\ti\lambda\leq\lambda-1$. Since
$(\sigma,S,\lambda)\in\Sigma_d$, we have $\sigma\neq 0$. Having fixed
the value of $\lambda$, we now run another induction procedure, this
one is downwards on the cardinality of~$S$. The base $|S|={n\choose
  d}$ is provided by Lemma~\ref{lm:easy}(3), since the left hand side
of \eqref{eq:est} is then equal to~$0$.  We now make the induction
step in~$|S|$.

Let us choose a~$(d-1)$-simplex $e\in\supp\sigma$. If $\sigma\in
B_{d-1}(\da;\zz)$, then there must exist a $d$-simplex $\tau\in\da(d)$
such that $e\in\bo\tau$. Let $\Omega$ denote the set of all
$d$-simplices $\tau\in\da_n(d)$ such that $e\in\bo\tau$. Clearly, we
have $|\Omega|=n-d$. We represent $\Omega$ as a~disjoint union
$\Omega=A\cup B\cup C$, where the sets $A$, $B$, and $C$ are defined
as follows:
\[A:=\{\tau\in\Omega\sm S\,\big{|}\,\,|\supp(\sigma+\bo\tau)|>(\lambda-1)(d+1)\},\]
\[B:=\{\tau\in\Omega\sm S\,\big{|}\,\,|\supp(\sigma+\bo\tau)|\leq(\lambda-1)(d+1)\},\]
\[C:=\Omega\cap S.\]

Since some simplex from $A\cup B$ must be picked in $\da$ we have the 
inequality
\begin{equation} \label{eq:in1}
\rho(\sigma,S,\lambda)\leq\sum_{\tau\in A\cup B}\prob(\tau)
\rho(\sigma+\bo\tau,S\cup\{\tau\},\lambda_\tau),
\end{equation}
where for each $\tau$ the value $\lambda_\tau$ is chosen so that
$(\sigma+\bo\tau,S\cup\{\tau\},\lambda_\tau)\in\Sigma_d$. In fact, we
shall see shortly that one can always choose $\lambda_\tau$ to be
$\lambda$ or $\lambda-1$. Substituting $p$ for $\prob(\tau)$, breaking
the sum on the right hand side of \eqref{eq:in1} into two, and using
the fact that $(\sigma+\bo\tau,S\cup\{\tau\},\lambda)\in\Sigma_d$ for
all $\tau\in A$, we obtain
\begin{equation} \label{eq:in2}
\rho(\sigma,S,\lambda)\leq 
p\sum_{\tau\in A}\rho(\sigma+\bo\tau,S\cup\{\tau\},\lambda)+
p\sum_{\tau\in B}\rho(\sigma+\bo\tau,S\cup\{\tau\},\lambda_\tau).
\end{equation}
Let $\alpha$ denote the first summand, and let $\beta$ denote the
second summand on the right hand side of~\eqref{eq:in2}.
We shall estimate these terms separately. 

First, since $|S\cup\{\tau\}|>|S|$, by the induction assumption (on
$|S|$) we have
\begin{multline} \label{eq:in3}
\alpha\leq p\,|A|\, c(d,\lambda)\, p^\lambda/(1-w)^\lambda<
p\, n\, c(d,\lambda)\, p^\lambda/(1-w)^\lambda=\\
=w\, c(d,\lambda)\, p^\lambda/(1-w)^\lambda.
\end{multline}

Let us next consider the summand $\beta$. To start with, if $\tau\in
B$, then $\supp\bo\tau$ contains at least one simplex from
$\supp\sigma$ other than $e$, and it is uniquely determined by that
simplex (together with~$e$). It follows that $|B|\leq|\supp\sigma|-1$.
Assume now that $\tau\in B$. In that case we have
\begin{multline} \label{eq:in4}
(\lambda-1)(d+1)\geq|\supp(\sigma+\bo\tau)|\geq|\supp\sigma|-|\supp\bo\tau|=
\\=|\supp\sigma|-(d+1)>(\lambda-2)(d+1),
\end{multline}
implying that $(\sigma+\bo\tau,S\cup\{\tau\},\lambda-1)\in\Sigma_d$,
and that $|\supp\sigma|\leq(\lambda+1)(d+1)$. Hence, by the induction
assumption (on $\lambda$) we have the estimate
\begin{multline}\label{eq:in5}
\beta\leq p\,|B|\, c(d,\lambda-1)\, p^{\lambda-1}/(1-w)^{\lambda-1}
\leq\\ \leq(|\supp\sigma|-1)\, c(d,\lambda-1)\, 
p^\lambda/(1-w)^{\lambda-1}<\\
<\lambda\,(d+1)\, c(d,\lambda-1)\, 
p^{\lambda}/(1-w)^{\lambda-1}.
\end{multline}

Substituting the estimates from \eqref{eq:in3} and \eqref{eq:in5} into
\eqref{eq:in2} we obtain
\begin{equation} \label{eq:in6}
\rho(\sigma,S,\lambda)< (w\, c(d,\lambda)+\lambda\,(d+1)\,
c(d,\lambda-1)\, (1-w))\, p^{\lambda}/(1-w)^{\lambda}.
\end{equation}
This yields the desired inequality \eqref{eq:est} for the constant
$c(d,\lambda)$ recursively defined by the equation
\[
c(d,\lambda):=w\, c(d,\lambda)+\lambda\,(d+1)\, c(d,\lambda-1)
\,(1-w),
\]
that is
\[
c(d,\lambda):=(d+1)\,\lambda\, c(d,\lambda-1).
\] 
Since $c(d,0)=1$, we arrive at
\begin{equation}\label{eq:in7}
c(d,\lambda)=(d+1)^{\lambda}\lambda!,
\end{equation}
which finishes the proof of the lemma.
\qed
\vspace{5pt}

We are now ready to proceed with the proof of the first part of 
our main theorem.

\vspace{5pt}

{\noindent{\bf Proof of Theorem~\ref{thm:main}(1). }} 

\nin Let us first settle the case $d=1$, as it can be done completely
explicitly, without referring to Lemma~\ref{star}. Clearly, for the
first Betti number of $\da\in Y_{n,p,1}$ to be nontrivial the graph
$\da$ must contain cycles. For $l=3,\dots,n$, let $z_l$ denote the
number of the $l$-cycles in a~complete graph on $n$ vertices. Then we
have
\begin{equation}\label{eq:d1}
\prob(\beta_1(Y_{n,p,1};\zz))\leq\sum_{\text{cycles\,\,} c}\prob(c)
=\sum_{l=3}^n z_l\, p^l.
\end{equation}
Substituting $z_l=\frac{1}{2}{n\choose l}(l-1)!$ into \eqref{eq:d1}
we obtain
\begin{multline}\label{eq:d2}
\prob(\beta_1(Y_{n,p,1};\zz))\leq\sum_{l=3}^n{n\choose l}(l-1)!
\, p^l/2=\\
=\sum_{l=3}^n \frac{n(n-1)\dots(n-l+1)}{2l}\, p^l
<\sum_{l=3}^n n^l\, p^l=w^3+\dots+w^n=\\
=w^3(1+w+\dots+w^{n-3})<\frac{w^3}{1-w}.
\end{multline}
In particular, 
\[\lni\prob(\beta_1(Y_{n,p,1};\zz))\leq\lni\frac{w(n)^3}{1-w(n)}=0.\]

For the rest of the proof we assume that $d\geq 2$. For an arbitrary
$d$-simplex $t$, let $A_t$ denote the event in $\ynpd$ that the chosen
complex $\da$ has a~nontrivial homology cycle which has
a~representative $\tau$ satisfying $t\in\supp\tau$. Let $t_0$ denote
the $d$-simplex with vertices $\{1,\dots,d+1\}$. Clearly, due to
symmetry, $\prob(A_t)=\prob(A_{t_0})$, for all $t\in\da_n(d)$, and so
we have
\begin{multline} \label{eq:t1}
\prob(\beta_d(\ynpd;\zz)>0)\leq\sum_{t\in\da_n(d)}\prob(t)\,\prob(A_t)=\\
=p\,{n\choose d+1}\,\prob(A_{t_0}) <
p\,n^{d+1}\,\prob(A_{t_0})=w\,n^d\,\prob(A_{t_0}).
\end{multline}

We shall next estimate $\prob(A_{t_0})$. As a~precursor of the general
argument we consider the case $d=2$. In this case $t_0$ is the triangle
with vertex set $\{1,2,3\}$. Let $e$ denote the edge with vertices $1$
and $2$. In order for the event $A_{t_0}$ to occur, we must pick some
triangle $s_i$ with the vertex set $\{1,2,i\}$, where $i=4,\dots,n$.
Hence we have the inequality
\begin{equation}\label{eq:t2}
\prob(A_{t_0})\leq\sum_{i=4}^{n}\prob(s_i)\,
\rho(\bo(t_0+s_i),\{s_i\},2).
\end{equation}
Since $|\supp(\bo(t_0+s_i))|=4>1\cdot 3=(\lambda-1)(d+1)$, by
Lemma~\ref{star} we have 
\[\rho(\bo(t_0+s_i),\{s_i\},2)\leq\frac{3^2\cdot 2!\cdot p^2}{(1-w)^2}
=\frac{18p^2}{(1-w)^2}.\] Combining this with \eqref{eq:t1}
and~\eqref{eq:t2}, and the fact that $\prob(s_i)=p$, we obtain
\[\prob(\beta_d(\ynptwo;\zz)>0)\leq w\, n^2\sum_{i=4}^{n}p\frac{18p^2}{(1-w)^2}
<\frac{18 w n^3 p^3}{(1-w)^2}=\frac{18 w^4}{(1-w)^2},\]
hence $\lni\prob(\beta_d(\ynptwo;\zz)>0)=0$ if $\lni w(n)=0$.

Let us now consider the case $d\geq 3$. The argument is along the same
lines as for $d=2$, but with more technical estimates, as it does not
suffice anymore to just add one $d$-simplex to $t_0$. Let
$e_1,\dots,e_{d+1}$ denote the $(d-1)$-dimensional faces of $t_0$
taken in an~arbitrary order. For the event $A_{t_0}$ to occur, for
each $i\in[d+1]$, we must pick at least one $d$-simplex different from
$t_0$ whose boundary contains~$e_i$.

Assume $T=\{t_1,\dots,t_{d+1}\}$ is such a~collection of
$d$-simplices, i.e., for all $i\in[d+1]$ we have
$t_i\in\da(d)\sm\{t_0\}$, and $e_i\in\supp(\bo t_i)$. For any
$i,j\in[d+1]$, $i\neq j$, we have $t_i\neq t_j$, since the only
$d$-simplex whose boundary contains both $e_i$ and $e_j$ is~$t_0$. We
consider the $d$-chain $\tau:=\sum_{i=0}^{d+1}t_i$.

Every $d$-simplex $t_i$ has a~unique vertex $v_i$ which does not
belong to $e_i$. We define a~set partition
$\pi=\pi_1\cup\dots\cup\pi_m$ on $T$ by putting $t_i$ and $t_j$ to the
same block if $v_i=v_j$.

\nin {\bf Claim.} {\it We have}
\begin{equation}\label{eq:t3}
|\supp(\bo\tau)|>(m-2)(d+1).
\end{equation}
\nin {\bf Proof of the Claim.} Clearly, $\supp(\bo\tau)$ consists of all 
the elements in $\bigcup_{i=0}^{d+1}\supp(\bo t_i)$ which belong to
the odd number of sets in that union. By construction, all the
elements of $\supp(\bo t_0)$ belong to exactly one other set in that
union, so all these cancel out. 

Potentially, we have $d(d+1)$ remaining elements. There will be no
cancellation between the elements of $\supp(\bo t_i)$ and $\supp(\bo
t_j)$ if $t_i$ and $t_j$ belong to different blocks in $\pi$. If they
belong to the same block, then there is exactly one cancellation,
namely of the $(d-1)$-simplices $\{v\}\cup(t_i\cap t_j)$, where $v$ is
the vertex corresponding to the block of $\pi$ containing $t_i$ and
$t_j$. Furthermore, all these cancellations are disjoint from each
other, since there are precisely two $(d-1)$-simplices in $\bo t_0$
containing $t_i\cap t_j$. We conclude that 
\begin{multline}\label{eq:tc1}
|\supp(\bo\tau)|=d(d+1)-\sum_{i=1}^m 2{|\pi_i|\choose 2}=
d(d+1)-\sum_{i=1}^m |\pi_i|(|\pi_i|-1)=\\
=d(d+1)+\sum_{i=1}^m|\pi_i|-
\sum_{i=1}^m |\pi_i|^2=(d+1)^2-\sum_{i=1}^m|\pi_i|^2.
\end{multline}
Since the sum $\sum_{i=1}^m|\pi_i|$ is fixed and all the terms in that
sum are positive integers, the maximum of $\sum_{i=1}^m|\pi_i|^2$ is
achieved by the values $|\pi_1|=\dots=|\pi_{m-1}|=1$,
$|\pi_m|=d+1-(m-1)$. Hence~\eqref{eq:tc1} yields
\begin{multline}
|\supp(\bo\tau)|\geq(d+1)^2-(m-1)-(d+1-(m-1))^2=\\
=(d+1)^2-(m-1)-(d+1)^2+2(d+1)(m-1)-(m-1)^2=\\
=(m-1)(2d+2-m)\geq(m-1)(2d+2-(d+1))>(m-2)(d+1),
\end{multline}
hereby proving \eqref{eq:t3}.
\qed

\vspace{5pt}

Since for $A_{t_0}$ to occur some constellation $T$ must be present in
our complex, we have an estimate
\begin{equation}\label{eq:t4}
\prob(A_{t_0})\leq\sum_\pi (n-d-1)^m\, p^{d+1}\,
\rho(\bo\tau,\supp\tau,m-1),
\end{equation}
where the sum is taken over all partitions
$\pi=\pi_1\cup\dots\cup\pi_m$, the factor $(n-d-1)^m$ records choosing
the $m$ vertices corresponding to the blocks of $\pi$, the factor
$p^{d+1}$ records the probability of choosing the set $T$, which is
uniquely determined by the choice of these vertices, and the term
$\rho(\bo\tau,\supp\tau,m-1)$ is well-defined by the claim which we
just proved, and the fact that $m\geq 1$. Using the
inequality~\eqref{eq:est}, we arrive at
\begin{multline}
\prob(A_{t_0})\leq\sum_\pi (n-d-1)^m p^{d+1}
(d+1)^{m-1}(m-1)!\, p^{m-1}/(1-w)^{m-1}<\\
<\frac{(d+1)^d\, d!}{(1-w)^d}\,p^d\sum_\pi n^m \, p^m
=\frac{(d+1)^d\, d!}{(1-w)^d}\, p^d\sum_\pi w^m
<\frac{(d+1)^d\, d!\,\npart(d+1)}{(1-w)^d}\,\,w\, p^d,
\end{multline}
where $\npart(d+1)$ denotes the number of set partitions of the set
$[d+1]$. Combining with \eqref{eq:t1}, end setting $c:=(d+1)^d\,
d!\,\npart(d+1)$, this yields
\[\prob(\beta_d(\ynpd;\zz)>0)<w\, n^d\, \frac{c}{(1-w)^d}\, w\, p^d=
c\,\frac{w^{d+2}}{(1-w)^d},\] We conclude that
$\lni\prob(\beta_d(\ynpd;\zz)>0)=0$ if $\lni w(n)=0$, also for all $d\geq
3$.
\qed

\section{Proof of the second part of Theorem~\ref{thm:main}}

Before we present the proof of the second part of
Theorem~\ref{thm:main}, we need to recall some standard tools of
combinatorial probability from~\cite{AS}. More specifically, a~certain
application of Chebyshev inequality has come to be known as the Second
Moment Method. We need the symmetric version of that method, which we
now proceed to describe.

Consider an infinite sequence of probability spaces $\cp^n$, where $n$
is a~natural number. Let us fix $n$ for now, and assume that we have
random events $A_1^n,\dots,A_m^n$ in $\cp^n$. For $i\in[m]$, let
$X_i^n$ denote indicator random variable of $A_i^n$, and set
$X^n=\sum_{i=1}^m X_i^n$. Assume furthermore, that the events $A_i^n$
are {\it symmetric} in the following sense: for every $i,j\in[m]$,
$i\neq j$, there exists an~automorphism of the underlying probability
space sending event $A_i^n$ to event $A_j^n$. An example of such
symmetric events in $\ynpd$ can be found by setting
$m:={n\choose{d+1}}$, indexing the $d$-simplices with the set $[m]$,
and letting $A_i^n$ denote the event that the $d$-simplex indexed with
$i$ lies in the chosen simplicial complex.

For distinct indices $i,j\in [m]$, we write $i\sim j$, in case the
events $A_i^n$ and $A_j^n$ are not independent. Furthermore, we set
\begin{equation}\label{eq:smm1}
\xi:=\sum_{i\sim j}\prob(A_i^n\wedge A_j^n).
\end{equation}
We mention explicitly that the sum in \eqref{eq:smm1} is taken over
all ordered pairs $(i,j)$, that is if the summand $\prob(A_i^n\wedge
A_j^n)$ occurs in the sum, then the summand $\prob(A_j^n\wedge A_i^n)$
occurs in the sum as well, since $i\sim j$ if and only if $j\sim
i$. Since for all $i,j\in[m]$ we have $\prob(A_i^n\wedge
A_j^n)=\prob(A_i^n)\prob(A_j^n|A_i^n)$, the equation~\eqref{eq:smm1}
now yields
\begin{equation}\label{eq:smm2}
\xi=\sum_{i=1}^m\prob(A_i^n)\sum_{j:i\sim j}\prob(A_j^n|A_i^n).
\end{equation}
We set 
\begin{equation}\label{eq:da*}
\xi^*:=\sum_{j:i\sim j}\prob(A_j^n|A_i^n),
\end{equation} which is
well-defined, since that sum does not depend on the choice of~$i$.

\noindent The following result can be found in \cite{AS}.
\begin{lm} \label{lm:smm} (\cite[Corollary 3.5]{AS})
With the notations above, if $\lni E(X^n)=\infty$ and
$\xi^*=o(E(X^n))$, then 
\begin{equation}\label{eq:smm3}
\lni\prob(X^n>0)=1,
\end{equation} 
and, furthermore,
\begin{equation} \label{eq:smm4}
\lni X^n/E(X^n)=1.
\end{equation}
\end{lm}

\noindent We now have all the necessary tools to proceed with the
proof of Theorem~\ref{thm:main}(2).

\vskip5pt

{\noindent{\bf Proof of Theorem~\ref{thm:main}(2).}} Our argument is
a~direct application of the second moment method. For fixed $d$ and
$p$, we set $\cp^n:=\ynpd$. We let
$\{\tau_1^n,\dots,\tau_{n\choose{d+2}}^n\}$ be the set of all
$(d+1)$-simplices of $\da_n$. For all $i=1,\dots,{n\choose{d+2}}$, let
$A_i^n$ denote the event that $\da\in\ynpd$ contains the boundary of
$\tau_i^n$, i.e., $\da(d+1)\supseteq\supp(\bo\tau_i^n)$.

As above, let $X_i^n$ denote the corresponding indicator random
variables, and set again
$X^n:=X_1^n+\dots+X_{n\choose{d+2}}^n$. Clearly,
$E(X_i^n)=\prob(A_i^n)=p^{d+2}$, for all $i$. Hence
\begin{multline*}
E(X^n)=\sum_{i=1}^{n\choose{d+2}}E(X_i^n)={n\choose{d+2}}p^{d+2}>
\frac{(n-d-1)^{d+2}\,p^{d+2}}{(d+2)!}=\\
=\frac{1}{(d+2)!}\left(1-\frac{d+1}{n}\right)^{d+2} w^{d+2},
\end{multline*}
and so we see that $E(X^n)=\Omega(w^{d+2})$, and, in particular, $\lni
E(X^n)=\infty$.

Furthermore, we have $i\sim j$ if and only if the $(d+1)$-simplices
$\tau_i$ and $\tau_j$ share precisely one boundary $d$-simplex. Thus,
in this case the dependency graph has $n\choose{d+2}$ vertices and is
regular of valency $(d+2)(n-d-1)$.

Given $i,j\in\left\{1,\dots,{n\choose{d+2}}\right\}$, such that $i\sim
j$, we get $\prob(A_j|A_i)=p^{d+1}$, since
$\left|\supp(\bo\tau_j)\sm\supp(\bo\tau_i)\right|=d+1$. Plugging this
data into the definition \eqref{eq:da*}, we get
\[\xi^*=\sum_{j:j\sim i}p^{d+1}=(d+2)(n-d-1)p^{d+1}.\]
Since
\[E(X^n)={n\choose{d+2}}p^{d+2}>\dfrac{n(n-d-1)^{d+1}}{(d+2)!}p^{d+2},\]
we get
\begin{equation}\label{eq:smm5}
\xi^*/E(X^n)<\frac{(d+2)(d+2)!}{np(n-d-1)^d}=\frac{(d+2)(d+2)!}{w(n-d-1)^d}.
\end{equation}
Since we assumed that $\lni w(n)=\infty$, the inequality
\eqref{eq:smm5} yields $\lni\xi^*/E(X^n)=0$, i.e., $\xi^*=o(E(X^n))$.
It then follows from Lemma~\ref{lm:smm} that $\lni\prob(X^n>0)=1$. 

Since $X^n>0$ implies that $\beta_d(\da;\zz)>0$, we get
$\prob(\beta_d(\da;\zz)>0)\geq\prob(X^n>0)$, hence
$\lni\prob(\beta_d(\ynpd;\zz)>0)=1$.  \qed


\vskip5pt

\section{Threshold probability for top homology group with coefficients
in an arbitrary abelian group} 
\label{sect:5}

In this short final section we shall indicate how to adjust our proofs
in order to deal with the case of homology with coefficients in
an~arbitrary abelian group. The exact statement which we get is the
following.

\begin{thm}
Assume that $p=p(n)=w(n)/n$, $d\geq 1$, and $F$ is an~arbitrary
nontrivial abelian group, then we have
\begin{enumerate}
\item if $\lni w(n)=0$, then $\lni\prob(H_d(\ynpd;F)=0)=1$;
\item if $\lni w(n)=\infty$, then $\lni\prob(H_d(\ynpd;F)\neq 0)=1$.
\end{enumerate}
\end{thm}

To start with, we need a~new piece of notations: for a~subset
$T\subseteq\da_n(d)$ we let $r(T)$ denote the number of
$(d-1)$-simplices $\sigma$ for which there exists a~unique $\tau\in T$
such that $\sigma\in\supp\bo\tau$. One may intuitively think of such
$(d-1)$-simplices as the ``rim'' of the set~$T$.

Next, the set $\Sigma_d$ should be replaced with
$\wti\Sigma_d\subseteq 2^{\da_n(d)}\times 2^{\da_n(d)} \times\cz_{\geq
  0}$ defined by the following: $(T,S,\lambda)\in\wti\Sigma_d$ if and
only if $r(T)>(\lambda-1)(d+1)$.

Accordingly, Definition~\ref{df:sigma} should be altered. For
$(T,S,\lambda)\in\wti\Sigma_d$, we define $\ti\rho(T,S,\lambda)$ to be
the probability that $\da\in\ynpd$ satisfies the following two
conditions:
\begin{enumerate}
\item[(1)] $\da(d)\cap S=\emptyset$;
\item[(2)] there exists $\sigma\in Z_d(\da_n)$ such that
  $T\subseteq\supp\sigma\subseteq T\cup\da(d)$.
\end{enumerate}

With this notations the inequality \eqref{eq:est} gets replaced with
\begin{equation} \label{eq:est2}
\ti\rho(T,S,\lambda)\leq 
c(d,\lambda)\, p^{\lambda}/(1-w)^{\lambda},
\end{equation}
with the proof holding almost verbatim. Essentially $|\supp\sigma|$
should be replaced with $r(\supp\sigma)$, and
$\rho(\bo\tau,S,\lambda)$ should be replaced with
$\ti\rho(\supp\tau,S,\lambda)$. For example, in the definition of $A$
and $B$ the expression $|\supp(\sigma+\bo\tau)|$ should be replaced
with $r(T\cup\{\tau\})$, the inequality~\eqref{eq:in1} becomes
\begin{equation} \label{eq:nin1}
\ti\rho(T,S,\lambda)\leq\sum_{\tau\in A\cup B}\prob(\tau)
\ti\rho(T\cup\{\tau\},S\cup\{\tau\},\lambda_\tau),
\end{equation}
and the chain of inequalities~\eqref{eq:in4} becomes 
\begin{multline} \label{eq:nin4}
(\lambda-1)(d+1)\geq r(T\cup\{\tau\})\geq r(T)-|\supp\bo\tau|=
\\=r(T)-(d+1)>(\lambda-2)(d+1),
\end{multline}
 
Also the proof of Theorem~\ref{thm:main}(1) holds verbatim with
similar changes. For example, the inequality \eqref{eq:t2} becomes
\begin{equation}\label{eq:nt2}
\prob(A_{t_0})\leq\sum_{i=4}^{n}\prob(s_i)\,
\ti\rho(\{t_0,s_i\},\{s_i\},2),
\end{equation}
the inequality~\eqref{eq:t3} in the claim becomes 
\begin{equation}\label{eq:nt3}
r(\{t_0,\dots,t_{d+1}\})>(m-2)(d+1),
\end{equation}
and the inequality~\eqref{eq:t4} becomes 
\begin{equation}\label{eq:nt4}
\prob(A_{t_0})\leq\sum_\pi (n-d-1)^m\, p^{d+1}\,
\ti\rho(\{t_0,\dots,t_{d+1}\},\{t_0,\dots,t_{d+1}\},m-1).
\end{equation}

Finally, the proof of Theorem~\ref{thm:main}(2) holds without any
changes at all since the presented $\zz$-cycles $\bo\tau_i^n$ are in
fact cycles for arbitrary coefficients.

\vskip5pt

\nin {\bf Acknowledgments.} The author would like to thank the
University of Bremen and the Banff International Research Station for
supporting this research.

\end{document}